\documentclass[reqno]{amsart}
\usepackage{fancyhdr}
\usepackage[dvips]{color}
\usepackage{mathrsfs}
\usepackage{amsfonts}
\usepackage{graphicx}
\usepackage{amscd}
\usepackage{amsmath}
\usepackage{amssymb}
\usepackage{latexsym}
\usepackage[all]{xy}
\usepackage[utf8]{inputenc}
\usepackage[T1]{fontenc}

\theoremstyle{plain}

\newtheorem{lemma}{\bf Lemma}

\newtheorem{proposition}{\bf Proposition}
\newtheorem{remark}{Remark}

\newtheorem{teo}{\bf Theorem}
\newtheorem{teo*}{Theorem}

\newcommand{\cqd}{{{\hfill q.e.d.}}\vspace{0.2cm}}

\numberwithin{equation}{section}
\def\CP{{\mathbb{CP}}}
\def\S{\mathbb{S}}

\def\C{\mathord{\mathbb C}}
\def\R{\mathord{\mathbb R}}

\def\g{\mathord{\frak{g}}}
\def\p{\mathord{\frak{p}}}

\def\Hom{\mathop{\rm Hom}\nolimits}
\def\10{_{(1,0)}}
\def\01{_{(0,1)}}

\def\2{\frac{1}{2}}
\def\Im{\mathop{\rm Im\,}\nolimits}
\def\Re{\mathop{\rm Re\,}\nolimits}
\def\.{\cdot}
\def\o{\circ}
\def\x{\times}
\def\Span{{\rm span}}

\def\beq{\begin{equation}}
\def\eeq{\end{equation}}
\def\bea{\begin{eqnarray}}
\def\eea{\end{eqnarray}}

\begin{document}

\title[Equivariant embeddings of kählerian symmetric spaces]{Equivariant embeddings of kählerian symmetric spaces}

\author{J.-H. Eschenburg}
\address{(J.-H. Eschenburg) Institut für Mathematik der
Universitäit Augsburg,
D-$86135$ Augsburg,
Germany}
\email{jost-hinrich.eschenburg@math.uni-augsburg.de}
\urladdr{https://myweb.rz.uni-augsburg.de/$\sim$eschenbu/}

\author{K.K. Santos}
\address{(K.K.Santos) Departamento de Matemática-CCT-UFRR, 69310-000, Boa Vista-RR-BR}
\email{kellykarina.kk@gmail.com}
\urladdr{http://www.ufrr.br}

\author{R. Tribuzy}
\address{(R. Tribuzy) Departamento de Matem\'atica-ICE-UFAM, 69080-900, Manaus-AM-BR}
\email{rtribuzy@yahoo.com.br}
\urladdr{http://www.ufam.edu.br}




\begin{abstract}
In this article we investigate some properties of 
equivariant embeddings of a symmetric Kählerian manifold. Motivated by a theorem of
Cartan and Wallach on equivariant embeddings of symmetric spaces we characterize 
these embeddings in the special case of $\mathbb{CP}^n$. 
Further, we verify that if a equivariant embedding has parallel plurimean curvature then it is the extrinsically symmetric one. 
\end{abstract}

\maketitle

\section{Introduction and statement of results}

The subject of the article is the theory of equivariant immersions of symmetric spaces, mainly Kähler symmetric spaces. A motivation for this investigation was the hope to find equivariant Kähler immersions with $\nabla^{\perp} \alpha^{(1,1)} = 0,$ so called ppmc (parallel pluri-mean curvature) immersions. 

The pluri-mean curvature $\alpha^{(1,1)}$ is the generalization of mean curvature in higher dimensions (\cite{Burs}). When $\alpha^{(1,1)}$ vanishes identically the immersion is called  $(1,1)$-geodesic or pluriminimal (\cite{Eells},\cite{Dajc1},\cite{Dajc2}). If  $\alpha^{(2,0)}$ vanishes identically, it is called $(2,0)$-geodesic which turns out to be the standard embedding of a Kähler symmetric space. 
Such immersions are obviously ppmc (\cite{Fe2}, \cite{EH Q}).

Though the hope for finding new ppmc immersions turned out not to be fulfilled, a new characterization theorem was found.

Let $P$ be a Riemannian manifold with a group $G$ of isometries acting transitively, a Riemannian homogeneous space. Denoting $K = G_p = \{g \in G : gp = p\}$ the isotropy subgroup (stabilizer) of some $p \in P,$ we may identify $P$ with the coset space $G/K$ by means of the $G$-equivariant diffeomorphism $P \ni gp \mapsto gK \in G/K.$ An immersion $\Phi : P \rightarrow V$ into some euclidean vector space V is called equivariant if there is a representation $\rho : G \rightarrow O(V )$ such that for all $g \in G$ and $p \in P$ we have
\begin{equation*}
\Phi(gp) = \rho(g)\Phi(p).
\end{equation*}
\noindent
In other words, $\Phi(P)$ is an orbit of the representation $\rho$, 
a very special orbit in fact: The
isotropy group $K = G_p$ fixes the vector $v_o = \Phi(p) \in V .$ A representation $\rho : G \rightarrow O(V )$ is called of class one for some subgroup $K \subset G$ if there is a nonzero fixed vector $v_o \in V$ for $K.$ Thus our representation
related to an equivariant immersion of $P = G/K$ is class-one for $K.$ Vice versa, if $K \subset G$ is a maximal subgroup, any representation of class-one for $K$ defines an equivariant {\em embedding} of $P$. We will restrict our attention to embeddings.

An equivariant embedding $\Phi : P \to V$ is called full if its image spans $V$. 
Another full equivariant embeddings $\tilde\Phi : P \to \tilde V$ is equivalent to $\Phi$
if there is a linear isometry $F : \tilde V \to V$ with $\tilde\Phi = \Phi \o F$.

A Riemannian manifold $P$ is called symmetric if at any point $p \in P$ there is an isometry $s_p$ of $P$ (the so called symmetry at $p$) fixing $p$ and reversing every geodesic through $p.$ A group of isometries $G$ containing the symmetry $s_p$ for every $p \in P$ acts transitively; hence $(P,G)$ is homogeneous. Further, $P$ is called Kähler symmetric if $s_p$ is part of a circle group $\mathbb{S}^1 \cong  S_p \subset K\subset G$ fixing $p$ and acting like the complex unit scalars on the tangent space $T_pP$ with some complex structure $J_p.$ In particular there is a distinguished square root of $s_p,$ an element $j_p \in S_p \subset G$ fixing $p$ with $j_p^2 = s_p.$

Kähler symmetric spaces $P = G/K$ enjoy a distinguished equivariant embedding 
$\Phi_o :P \rightarrow V_o = \g,$ the so called standard embedding which assigns to 
each $p \in P$ the infinitesimal generator $J_p$ of $S_p.$ This is extrinsically 
symmetric, that is its second fundamental form $\alpha$ is parallel, 
$\nabla^{\perp} \alpha = 0.$ Hence it is trivially ppmc: Not only the $(1,1)$-component
of $\nabla^\perp\alpha$ is zero, but the full derivative $\nabla^\perp\alpha$ vanishes.

However, as it was shown by \'Elie Cartan \cite{Ca} and Nolan Wallach \cite{Wallach}, there are infinitely many inequivalent equivariant embeddings for any symmetric space $P = G/K$. Let $\Phi : P \rightarrow V$ be such an embedding and $\rho : G \rightarrow O(V)$ the corresponding representation. As it turns out (see Prop.\ Lemma \ref{decompo}), it suffices to check the ppmc property only for a certain subset of equivariant embeddings, so called fundamental ones, which generate all the others.

Let $N = N_p\Phi(P) = d\Phi(T_pP)^{\perp}$ be the normal space at $p$. 
Then $N$ decomposes as $N^+ \oplus N^-$, the $\pm1$-eigenspaces of $\rho(s_p)|N$. 
If $P$ is Kähler symmetric, $N$ decomposes further into the eigenspaces of
$\rho(j_p)$ at $p.$ We denote $N^{++} \subset N^{+}$ the $1$-eigenspace of $\rho(j_p)|N.$ 
One important observation is (see the lemmas \ref{alpha} and \ref{alpha11lemma}):

\vspace{0.3cm}

\begin{proposition} \label{alphaprop}
The second fundamental form $\alpha$ takes values in $N^{+}$, 
and its $(1,1)$-part is the 
projection of $\alpha$ onto $N^{++}.$
\end{proposition}

\vspace{0.3cm}

We are mainly interested in the most basic Kähler symmetric space,
complex projective space $\CP^n$. 
In Prop.\ \ref{eqembprop} we construct all fundamental equivariant embeddings of $\CP^n$ 
in Euclidean space (even without using the theorem of Cartan-Wallach). 
Then we compute $\nabla^\perp\alpha^{(1,1)}$ for every such embedding, 
see Equation \eqref{nablaalpha11}. This yields our main result:

\begin{teo}\label{theoCPn}
None of the equivariant embeddings of $\mathbb{CP}^n$ is ppmc except the standard embedding.
\end{teo}


The present paper is a revised version of a part of the PhD thesis \cite{Santos} of the
second author.

\section{Eigenspace decomposition for a complex structure}  \label{eigen}

Let $T$ be a Euclidean vector space equipped with a complex structure,
that is an orthogonal linear map $J$ with $J^2 = -I$ 
(where $I$ denotes the identity map on $T$). Then $J$ has eigenvalues $\pm i$,
and $T^c := T\otimes\C$ decomposes orthogonally into the eigenspaces
$T\10 = \{v-iJv: v\in T\}$ (eigenvalue $i$) and $T\01 = \overline{T\10}$ (eigenvalue $-i$)
\beq	\label{type}
	T^c = T\10 \oplus T\01.
\eeq
Accordingly, for any $r\in\{0,1,2,\dots\}$ the tensor power 
$(T^c)^r 
= \bigotimes^rT^c$ splits as
$$
	(T^c)^r = \bigoplus_{p+q=r} (T\10)^p\otimes (T\01)^q.
$$
Let $\alpha : T^r \to N$ be a linear map into another vector space $N$. After extending it
to a complex linear map $\alpha : (T^c)^r \to V^c$, we let $\alpha^{(p,q)}$ be the
restriction to the subspace $(T\10)^p\otimes (T\01)^q$,
$$
	\alpha^{(p,q)} := \alpha|\left((T\10)^p\otimes (T\01)^q\right)
$$
for any $p,q$ with $p+q=r$.

We will apply this to the second fundamental form of an immersed Kähler manifold $P$.
Recall that a Kähler manifold is a Riemannian manifold $P$ with a parallel complex structure
$J$ on its tangent bundle $TP$. Then the decomposition \eqref{type} on tangent spaces
$T = T_pP$ is parallel. Now let $\Phi : P\to V$ be an isometric immersion into
Euclidean vector space $V = \R^n$ and $\alpha : TM\otimes TM \to NM$ its second fundamental
form. Then its complexification has the components $\alpha^{(2,0)}$, $\alpha^{(0,2)}$, and
$\alpha^{(1,1)}$. In particular, for all $v\in T$
$$
	\alpha^{(1,1)}(v,v) = \alpha(v\10,v\01) = \2\alpha(v-iJv,v+iJv)
	= \2(\alpha(v,v)+\alpha(Jv,Jv)).
$$
If $\dim P = 2$ and $|v|=1$, this is the mean curvature vector. Therefore $\alpha^{(1,1)}$
is called pluri-mean curvature.

\section{Equivariant embeddings of compact symmetric spaces}

Let $P = G/K$ be a compact symmetric space.
For any equivariant immersion $\Phi : P \rightarrow V$ , the vector space $V$ can
be viewed as a space of functions on $P,$ a subspace of $C^{\infty}(P).$ In fact,
we assign to each $v \in V$ the height function
\begin{equation*}
f_v : P \rightarrow \mathbb{R} : f_v(x) = \langle \Phi(x), v \rangle.
\end{equation*}
where $\langle\ ,\,\rangle$ denotes the Hermitian inner product $\langle v,w\rangle = v^*w$ which 
is antilinear on the first argument.

Then the $G$-representation $\rho$ on $V$ becomes part of the $G$-representation
on $C^{\infty}(P)$ given by transforming the argument of a function $f \in C^{\infty}(P).$ This defines a new function $gf \in C^{\infty}(P)$, for any $g \in G$ as follows:
\begin{equation} \label{gf}
(gf)(x):=f(g^{-1}x).
\end{equation}

To see compatibility of the representations on $V$ and $C^{\infty}(P)$ we observe
$$f_{\rho(g)v}(x)=\langle \Phi(x), \rho(g)v \rangle = \langle \rho(g)^{-1} \Phi(x),v \rangle= \langle \Phi(g^{-1}x),v \rangle = (gf_v)(x).
$$
Thus we have seen: 

\begin{proposition} \label{eqemb}
Every equivariant embedding of a symmetric space $P = G/K$ into Euclidean space 
is obtained from a finite-dimensional class-one $G$-submodule of $C^{\infty}(P)$.
This is the orbit $G.f = \{f\o g: g\in G\}$ for some function $f\in C^\infty(P)$
which is fixed precisely by the subgroup $K$ and has finite dimensional linear span.
\end{proposition}

A theorem of \'Elie Cartan and Wallach (\cite{Ca,Wallach}, cf.\ also \cite[App.\ B]{EHQ})
states that $C^{\infty}(P)$ decomposes completely into irreducible $G$-representations $V_k$ which are of class one with respect to $K$, and each such representation occurs precisely once in $C^{\infty}(P)$, up to equivalence. 
In other words, for any $k$ there is a ``fundamental'' 
$K$-invariant function $f_k$ on $P$ such that
$$
	V_k = \Span\,Gf_k.
$$
In view of Proposition \ref{decompo} below, it will be enough for our question (``ppmc'') 
to consider those equivariant immersions of $G/K$ which are orbits of $f_k$ for some $k$,
$$
	Gf_k=\{f_k \circ g^{-1}:g \in G\}.
$$

\section{Equivariant embeddings of complex projective space (1)}
Now we turn to our main example, complex projective $n$-space:
\beq	\label{CPn}
	P=\CP^n = \S^{2n+1}\!/\S^1 = \{[z]: z\in\S^{2n+1}\} = G/K
\eeq 
with $G = U_{n+1}$ and $K = U_n\x U_1$, where $[z] := \S^1\.z$. 
Our base point is $p = [e_{n+1}]$. One can generate 
the algebra $C^{\infty}(P)$ by certain real polynomials
on $\mathbb{C}^{n+1}$, viewed as polynomials in $z \in \mathbb{C}^{n+1}$ and its complex conjugate $\overline{z}$ (used as coordinates in place of the variables $\Re z$ and $\Im z$). These polynomials descend from $\mathbb{S}^{2n+1} \subset \mathbb{C}^{n+1}$ to $\mathbb{CP}^n$ if and only if they decompose into polynomials $f_k$ which are homogeneous of the same degree $k$ in $z$ and $\bar{z}$ (``bihomogeneous'' of degree $k$). Such polynomial $f_k$ is $K$-invariant for $K = U_n\times U_1$ if and only if it depends only on the last variable $z_{n+1}$,
\beq	\label{fk}
	f_k =  (e_{n+1}^* \overline{e_{n+1}^*})^k : 
	z \mapsto z_{n+1}^k \overline{z_{n+1}}^k = |z_{n+1}|^{2k}
\eeq
for $k\in\{0,1,2,\dots\}$. 
The $G$-orbit of $f_k : [z]\mapsto |(e_{n+1})^*z|^{2k}$ 
consists of the functions $gf_k = f_k\o g^{-1}$ with 
\beq \label{gfk}
	gf_k([z]) = f_k(g^*z)) = (e_{n+1}^*g^*z)^k = |(ge_{n+1})^*z|^{2k} = |v^*z|^{2k}
\eeq 
for $v = ge_{n+1}$, for any $g\in G = U_{n+1}$.

Since $f_k$ belongs to the space $V_k$ of $k$-bihomogeneous polynomials on $\C^{n+1}$
and since the positive real multples of $f_k$ are the only 
$K$-invariant polynomials in $V_k$, we have seen:

\begin{lemma} \label{Vk}
Let $V_k$ be the space of real bihomogeneous polynomials of bidegree $k$ on $\C^{n+1}$.
Then $V_k$ is the real span of the orbit $U_{n+1}f_k$ where $f_k$ is defined in {\rm \eqref{fk}}.
\end{lemma}

\begin{proposition} \label{eqembprop}
Every equivariant embedding of $\CP^n$ (up to equivalence) 
is an orbit $U_{n+1}.f= \{f\o g: g\in U_{n+1}\}$ where $f : \CP^n\to\R$ with
\beq \label{f=sumfk}
	f([z]) = \sum_k a_k |z_{n+1}|^{2k}
\eeq 
for all $z=(z_1,\dots,z_{n+1})\in\C^n$ with $|z|=1$, and finitely many $a_k\in\R$ are nonzero.
\end{proposition}

\proof
Let $G = U_{n+1}$ and $K = U_n\x U_1$.
By Prop.\ \ref{eqemb} we may assume that our equivariant embedding is the $G$-orbit 
of a smooth function $f$ on $\CP^n$ such that $f$ is precisely invariant under $K$ and
$Gf$ spans a finite dimensional $G$-submodule $V \subset C^\infty(\CP^n)$. We may
consider $f$ as an $\S^1$-invariant function on $\S^{2n+1}$.
Since $f$ can be uniformly approximated by $\S^1$-invariant 
polynomials on $\S^{2n+1}$, it is not perpendicular to $\sum_k V_k$ with respect to the
$L^2$-inner product. Thus $\pi_k(f) \neq 0$ for some $k$, where $\pi_k$ denotes the orthogonal
projection onto $V_k$.
But the submodules $V_k = \Span(Gf_k)$
are mutually inequivalent since they differ by dimensions. 
Thus $V = \Span(Gf) = \sum_k \Span(G.a_kf_k)$.
Since $V$ is finite dimenisonal, only finitely many $a_k$ can be nonzero
which proves \eqref{f=sumfk}.
\endproof

\section{Equivariant embeddings of a Kähler Symmetric Space}


Let $P$ be any compact symmetric space and $\Phi: P \rightarrow V$ an equivariant embedding. If $p \in P$ so we will use $T$ to denote $T_pP$ and $N$ to $N_pP$. Let $s$ be the symmetry in $p$ and $\rho(s)$ its extension to an orthognal space of $V$.
Let $E_\pm(s)$ the the $(\pm 1)$-eigenspaces of $\rho(s)$ and set
$$
	N^+ := E_+(s)\subset N, \ \ \ N_- := E_-(s)\cap N.
$$
Note that $E_+(s)\subset N$ since $s = -I$ on $T$. Then $N =  N^+ + N^-$.

\begin{lemma} \label{alpha}
The second fundamental form $\alpha$ of $\Phi$ takes values in $N^+$, and
$\nabla^\perp\alpha$ takes values in $N^-$.
\end{lemma}

\proof
In fact, for all $v,w\in T$ we have
\begin{equation*}
s \alpha(v,w)=\alpha (sv,sw)=\alpha (-v,-w) = \alpha (v,w),
\end{equation*}
that is, $\alpha (V,W) \in N^{+}$. Analogously, for all $u,v,w\in T$,
\begin{equation*}
s \nabla_u^{\perp} \alpha(v,w)= \nabla_{su}^{\perp} \alpha (sv,sw)= \nabla_{-u}^{\perp} \alpha (-v,-w) = - \nabla_u^{\perp} \alpha (v,w),
\end{equation*}
that is, $\nabla_u^{\perp} \alpha (v,w) \in N^{-}$. 
\endproof


\begin{remark}
If $P$ is a Kähler symmetric space and $\Phi: P \rightarrow V$ is an equivariant embedding so $\langle \alpha^{(2,0)},\alpha^{(2,0)} \rangle$ is holomorphic.
\end{remark}

\begin{dem}
In fact, let $U,V,W$ and $Z$ be $(1,0)$-vector fields with zero derivative at the base point $p$, and let $A$ be a vector of type $(0,1)$. Then 
$$\begin{array}{lll}
 A\langle \alpha(U,V), \alpha(W,Z) \rangle   & = & \langle \nabla_{A}^{\perp} \alpha(U,V), \alpha(W,Z) \rangle +
                                                        \langle  \alpha(U,V), \nabla_{A}^{\perp} \alpha(W,Z) \rangle \\
                                                   & = & \langle (\nabla_{A}^{\perp} \alpha)(U,V) + \alpha (\nabla_{A}U,V) + \alpha (U,\nabla_{A}V)\,,\, \alpha(W,Z) \rangle +\\
                                                   &   & \langle \alpha(U,V)\,,\, (\nabla_{A}^{\perp} \alpha)(W,Z) + \alpha (\nabla_{A}W,Z) + \alpha (W,\nabla_{A}Z) \rangle\\
                                                   & = & 0
 \end{array}$$
The last equality is a consequence of $\nabla_{A}U= \nabla_{A}V= \nabla_{A}W = \nabla_{A}Z=0$ and the fact that 
$\nabla ^{\perp} \alpha$ and
$\alpha$ take values in $N^+$ and $N^-$ (see Lemma \ref{alpha}) which are perpendicular. 
 \end{dem}

 \vspace{0.2cm}
 
Therefore if $\Phi:P \rightarrow V$ is an equivariant embedding so $\alpha \in N^+$ and $\nabla^{\perp} \alpha \in N^-$, in other words, they are orthogonal. On the other hand we have:

\begin{remark} 
If  $\Phi: P \rightarrow V$ is an embedding with $\alpha$ and $\nabla^{\perp} \alpha$ 
orthogonal so the \linebreak embedding is locally symmetric and if  $P$ is simply connected the embedding is \linebreak symmetric, $($\cite{Ziller}, pg 138$)$.
\end{remark}

\vspace{0.3cm}

\begin{dem}
By the Gauss equation 

$\begin{array}{lll}
 \langle R(X,Y)Z, W \rangle & = & \langle \alpha(X,Z),\alpha(Y,W) \rangle - \langle \alpha(Y,Z),\alpha(X,W) \rangle
 \end{array}$.
So,

$\begin{array}{lll}
\nabla \langle R(X,Y)Z, W \rangle & = & \langle \nabla^{\perp} \alpha(X,Z),\alpha(Y,W) \rangle + \langle \alpha(X,Z),\nabla^{\perp} \alpha(Y,W) \rangle \\
                            &   &- \langle \nabla^{\perp} \alpha(Y,Z),\alpha(X,W) \rangle - \langle \alpha(Y,Z),\nabla^{\perp} \alpha(X,W) \rangle \\
                            & = & 0 \;\;  (\mbox{because} \; \nabla^{\perp} \alpha \; \mbox{and} \; \alpha \; \mbox{are perpendicular}).
 \end{array}$
 
\smallskip\noindent
Hence the curvature is covariantly parallel, i.e.\ the imersion is locally 
symmetric. \cqd

\end{dem}

\vspace{0.2cm}

Now let be $P$ Kähler symmetric and $J$ the complex structure on $P$. So $J$ extends to an isometry on $V$ with $J^2 = s.$ Since $s = I$
on $N^+$, this space splits as $N^{++} + N^{+-}$  where $J = 1$ on $N^{++}$ and $J = -1$ on $N^{+-}.$

Let be $V_{(1,0)},W_{(1,0)}$ of type $(1,0)$ and $V_{(0,1)},W_{(0,1)}$ of type $(0,1)$, so


\beq	\label{alphaV}
\begin{matrix}
J \alpha (V_{(1,0)},W_{(0,1)})  & = & \alpha (J V_{(1,0)},J W_{(0,1)}) \cr
	& = & \alpha(iV_{(1,0)},-iW_{(0,1)}) & = & \alpha(V_{(1,0)},W_{(0,1)})\\
J \alpha (V_{(0,1)},W_{(1,0)}) & = & \alpha (J V_{(0,1)},J W_{(1,0)})\cr 
	& = & \alpha(-iV_{(0,1)},iW_{(1,0)}) & = & \alpha(V_{(0,1)},W_{(1,0)})\\
J \alpha (V_{(1,0)},W_{(1,0)})  & = & \alpha (J V_{(1,0)},J W_{(1,0)}) \cr
	 & = & \alpha(iV_{(1,0)},iW_{(1,0)})  &= & - \alpha(V_{(1,0)},W_{(1,0)})   \\
J \alpha (V_{(0,1)},W_{(0,1)}))& = & \alpha (J V_{(0,1)},J W_{(0,1)}))\cr
	 &= & \alpha(-iV_{(0,1)}),-iW_{(0,1)})   &= & - \alpha(V_{(0,1)}),W_{(0,1)})
\end{matrix}
\eeq

\vspace{.2cm}

Since $\alpha : T\otimes T\to N$ takes values in $N^+ = N^{++}\oplus N^{+-}$, we have seen: 

\begin{lemma} \label{alpha11lemma}
Let $\pi^{++} : N^+ \to N^{++}$ and $\pi^{+-} : N^+ \to N^{+-}$ be the orthogonal projections.
Then
\beq\begin{matrix} \label{alphapq}
			\alpha^{(1,1)} &=& \pi^{++}\o\alpha,	\cr
	 \alpha^{(2,0)}+\alpha^{(0,2)} &=& \pi^{+-}\o\alpha.	
\end{matrix}\eeq
\end{lemma}

\proof
$$\begin{array}{lll}
\alpha (X,Y)\!\!&=&\!\!\alpha (X_{(1,0)}+X_{(0,1)},Y_{(1,0)}+Y_{(0,1)}) \\
             \!\!&\buildrel (*)\over=&\!\!\alpha (X_{(1,0)},\!Y_{(1,0)})+ \alpha (X_{(0,1)},\!Y_{(0,1)})
             + \alpha (X_{(0,1)},\!Y_{(1,0)}) + \alpha (X_{(1,0)},\!Y_{(0,1)}) \\
             \!\!&\buildrel (**)\over=&\!\! \alpha^{(2,0)}(X,Y) + \alpha^{(0,2)}(X,Y)+\alpha^{(1,1)}(X,Y).
\end{array}$$
By \eqref{alphaV}, the  first two terms on the right hand side of $\buildrel (*)\over=$   belong to $N^{+-}$ 
while the last two terms belong to $N^{++}$. The equality $\buildrel (**)\over=$ shows that this is 
the splitting of $\alpha$ into $(\alpha^{(2,0)}+\alpha^{(0,2)})$ and $\alpha^{(1,1)}$. 
\endproof

We will see later for $P = \CP^n$ that no equivariant embedding in $V_k$, $k \geq 2$, 
has parallel pluri-mean curvature. To verify that the same holds for equivariant embeddings in sums of several $V_k$ 
we will use the following general statement:


\begin{proposition} \label{decompo}
Let $P$ be any Kähler manifold and $\Phi' : P \to V'$, $\Phi'' : P \to V''$ two
isometric immersions in Euclidean vector
spaces $V',V''$. Suppose that $\Phi'' = \Phi\oplus\Phi' :
P \to V := V'\oplus V''$ is a ppmc immersion. Then $\Phi',\Phi''$ are also ppmc.
\end{proposition}

\proof 
The normal space $N$ of $\Phi$ at some $p\in P$ splits
as the sum of the normal spaces $N' \oplus N''$
and the tangent antidiagonal space $N^o = \{v'\oplus -v'': v\in T\}$
where we let $v' = d\Phi'_p(v)$ and $v'' = d\Phi''_p(v)$. 
This is a parallel and orthogonal splitting. Let $\alpha$ be the second
fundamental form of $\Phi$. Since $\alpha^{(1,1]}$
is parallel, its three components in $N',N'',N^o$ must be parallel,
in particular the $N'$ and $N''$ components, which are the $(1,1)$-components
of $\alpha'$ and $\alpha''$, the second fundamental forms of $\Phi'$ and $\Phi''$. 
Therefore both $\Phi'$ and $\Phi''$ are ppmc immersions.
\endproof

\section{Equivariant embeddings of complex projective space (2)}

Continuing our study of the equivariant embeddings of
complex projective $n$-space $\CP^n = \C^{n+1}_*/\C_* = \S^{2n+1}/\S^1$
(where $\C^{n+1}_* = \C^{n+1}\setminus\{0\}$) we will now prove Theorem 1.

To save notation will slightly modify our model for $\CP^n$ 
by replacing the vector space $\C^{n+1}$
with its dual vector space $(\C^{n+1})^* = \Hom_{\C}(\C^{n+1},\C)$. 
The elements of $(\C^{n+1})^*$ are $\C$-linear functions $x : \C^{n+1}\to\C$,
that is $x = v^* : z \mapsto v^*z$ with $v,z\in\C^{n+1}$, and the elements $[x]\in
\CP^n$ are the $\S^1$-orbits of $x$, 
that is $[x] = \{\lambda x: \lambda\in\S^1\subset\C\}$.

Now our base point $x_o$ in $\CP^n$ is the function $e_{n+1}^*$, and $x := gx_o = v^*$
where $v = ge_{n+1}$ is the last column of the matrix $g$. Hence 
$$
	\Phi_k(x) = |x\bar x|^k
$$
for the function $\bar x : z\mapsto x(\bar z)$.

By Prop.\ \ref{eqembprop} and Proposition \ref{decompo} we only need to consider this embedding
$\Phi = \Phi_k$ for any $k$. It takes values in the vector space $V_k$ of all polynomials of bidegree $k$
on $\C^{n+1}$ where the representation $\rho_k$ of $G = U_{n+1}$ 
on the function space $V_k$ as in \eqref{gf}.

How does the symmetry\footnote 
	{Recall that $[x] = \S^1x \in \CP^n$. 
	We will write $s_x$ for $s_{[x]}$ and $j_x$ for $j_{[x}]$.}
$s_x$ and its square root $j_{x}$ act on $V_k$?
Let $x(t)$ be a great circle in $\S^{2n+1}$ starting at $x(0) =:x$ with inital 
derivative $v = \delta x \perp \C x$  where 
	$$\delta x := \left.\frac{dx(t)}{dt}\right|_{t=0}.$$ 
Since $x(t)$ is a great circle,
	$$\delta v = -x.$$
If $x$ is the base point $(e_{n+1})^*$, 
then $s_x = (-I,1) \in U_n\x U_1$ where $I\in U_n$ denotes the
unit matrix on $\C^n$, and $j_x = (iI,1)$. Thus $j_x$ (more precisely, $(dj_x)_x$)
is the complex structure $iI$ on $\C^n \cong (\C x)^\perp=T_{[x]}\CP^n$ 
while $j_x(x) = x$. The same holds for any $[x]\in\CP^n$ 
since $G = U_{n+1}$ perserves the complex structure on $\CP^n$.
Then $\delta\Phi_k x = (d\Phi_k)_x\delta x$ is a polynomial expression in $x,\bar x,v,\bar v$,
and applying $\rho_k(s_x)$ or $\rho_k(j_x)$ we obtain the same expressions where every $v$ and $\bar v$ is replaced by 
$-v$ and $-\bar v$ or by $iv$ and $-i\bar v$, respectively.
 
\vspace{0.2cm}

Let us calculate the second fundamental form of $\Phi_k$. Let $x,x(t),v$ be
as above.
$$\
\hat v := d\Phi_{[x]} v= \delta (x\bar x)^k = k(x\bar{x})^{k-1}(v \bar{x}+x \bar{v}).
$$

Moreover,

\vspace{0.2cm}

$$\delta \hat v= k(k-1)(x \bar{x})^{k-2}\delta (x \bar{x})(v \bar{x}+x \bar{v})+ k(x \bar{x})^{k-1} \delta
(v \bar{x}+x \bar{v}). $$

\vspace{0.2cm}

Using $\delta v = -x$ we have
\bea
\delta(x \bar{x})&=&v \bar{x}+ x \bar{v},\cr
\delta (v \bar{x}+x \bar{v})&=& 2(- x \bar{x}+ v \bar{v}).
\nonumber
\eea 

\vspace{0.2cm}
Hence

$$\begin{array}{lll}
\delta \hat v & = &  k(k-1)(x \bar{x})^{k-2}(v \bar{x}+x \bar{v})^2+2k(x \bar{x})^{k-1}(-x \bar{x}+v \bar{v})  \\
         & = & k(k-1)(x \bar{x})^{k-2}(\bar{x}\bar{x}vv+xx \bar{v}\bar{v}+2x\bar{x}v \bar{v})+2k(x \bar{x})^{k-1}(-x \bar{x}+v \bar{v})
 \end{array}$$

\vspace{0.2cm}

Furthermore 
$\alpha ^{( 1,1)}(v,v)= \left( \delta \hat v \right)^{N^{++}}$. Thus we have to project $\delta\hat v$ into $E_+(j_x)$, the
$(+1)$-eigenspace of $j_x$. This is automatically contained in the $+1$-eigenspace of $s_x$ and thus in the normal space $N$
since $s_x = -I$ on $T$. A term of $\delta\hat v$ is contained in $E_+(j_x)$ if it contains the same number of factors 
of type $v$ and $\bar v$, otherwise it lies in $E_-(j_x)$. Thus precisely the terms with $\bar x\bar x vv$ or $xx\bar v\bar v$
are projected to zero, and we obtain
\beq \label{alpha11}
\alpha^{( 1,1)}(v,v)= 2k(k-1)(x \bar{x})^{k-2}x \bar{x} v \bar{v}+
2k(x \bar{x})^{k-1}(-x \bar{x} + v \bar{v}).
\eeq

\vspace{0.4cm}

If $k=1$ the expressions for $\delta \hat v = \alpha(v,v)$ and $\alpha ^{(1,1)}(v,v)$ coincide, in fact, both expressions are $2(-x \bar{x}+v \bar{v})$. 

\vspace{0.4cm}

To calculate $\nabla_v \alpha^{(1,1)}(v,v)$ we have to differentiate the right hand side of \eqref{alpha11}, the normal vector 
$\xi =2k(k-1)(x \bar{x})^{k-2}x \bar{x} v \bar{v}+
2k(x \bar{x})^{k-1}(-x \bar{x} + v \bar{v})$, and take its component in $N^-$ afterwards: 

\vspace{0.2cm}

$$\begin{array}{lll}
\delta \xi & = & 2k(k-1)(k-2)(x \bar{x})^{k-3}(v\bar{x}+x \bar{v})(x \bar{x}v \bar{v})\,+\\
           &   & 2k(k-1)(x \bar{x})^{k-2}((v \bar{x}+x \bar{v})v \bar{v} + x \bar{x}(-x \bar{v}- v \bar{x}))\,+\\
           &   & 2k(k-1)(x \bar{x})^{k-2}(v \bar{x}+x \bar{v})(-x \bar{x}+v \bar{v}\,)\,+\\
           &   & 2k(x \bar{x})^{k-1}(-(v \bar{x}+x \bar{v})-x \bar{v}-v \bar{x})) \\
           &   & \\
	   & = & 2k(k-1)(k-2)(x \bar{x})^{k-3}(x \bar{x} \bar{x}vv \bar{v} + x x \bar{x}v \bar{v} \bar{v})\,+\\
           &   & 4k(k-1)(x \bar{x})^{k-2}(xv \bar{v} \bar{v}+ \bar{x} vv \bar{v} -xx \bar{x} \bar{v}-x \bar{x}\bar{x}v)\,+\\
           &   & 4k(x \bar{x})^{k-1}(-x \bar{v}-v \bar{x})
\end{array} $$

\vspace{0.2cm}

Recall that $N^-$ is the $(-1)$-eigenspace of $\rho_k(s_x)$ which contains all terms with an odd number of factors
of type $v$ or $\bar v$. Therefore the component in $N^-$ is:
\beq \label{nablaalpha11}
\begin{array}{lll}
\nabla_v \alpha^{(1,1)}(v,v) & = &  \left( \delta \xi \right)^{N^-} \\
                             & = & 2k(k-1)(k-2)(x \bar{x})^{k-3}(x \bar{x} \bar{x}vv \bar{v} + x x \bar{x}v \bar{v} \bar{v})\,+\\
                             &   & 4k(k-1)(x \bar{x})^{k-2}(xv \bar{v} \bar{v}+ \bar{x} vv \bar{v} -xx \bar{x} \bar{v}-x \bar{x}\bar{x}v)
\end{array}\eeq
%
%
The above expression is zero only for $k=1$, the extrinsically symmetric case. 
So we have 
proved our main result:

\vspace{0.2cm}

{\bf Theorem \ref{theoCPn}}
{\sl None of equivariant embeddings of $\mathbb{CP}^n$ is ppmc but the standard embedding.}

\vspace{0.5cm}

{\bf Acknowledgements:} The second author would like to thank the Department of Mathematics of Augsburg University, where part of this work was carriod out. The second author was been partially supported by CAPES-Brazil. The third author is partially supported by Conselho Nacional de Desenvolvimento Científico e Tecnológico (CNPq), of the Ministry of Science, Technology and Innovation of Brazil and Fundação de Amparo à Pesquisa do Estado do Amazonas (FAPEAM).



\end{document}